\documentclass[12pt]{article}
\usepackage{amsmath,amssymb,amscd,amsthm,latexsym}

\font\block=msbm10
\def\C{\hbox{\block\char'0103}}
\def\Z{\hbox{\block\char'0132}}

\font\gotic=eufm10
\def\g{\hbox{\gotic\char'0147}}
\def\l{\hbox{\gotic\char'0154}}
\def\s{\hbox{\gotic\char'0163}}
\def\p{\hbox{\gotic\char'0160}}
\def\o{\hbox{\gotic\char'0166}}

\def\h{\hbox{\gotic\char'0150}}

\font\new=eusm10

\def\P{\hbox{\new\char'0120}}

\font\yes=msam10
\def\Q{\hbox{\yes\char'0003}}

\begin{document}

\noindent
{\Large On cohomology of the Lie superalgebra $D(2, 1 ; \alpha)$}

\vskip 0.5in
{\large Elena Poletaeva}
\vskip 0.1in

{\it Department of Mathematics,}

{\it University of Texas-Pan American,}

{\it Edinburg, TX 78539}

{\it Electronic mail:} elenap$@$utpa.edu

\vskip 0.4in

{\footnotesize \noindent {\bf Abstract.}
We describe the infinitesimal deformations of the standard embedding of the Lie superalgebra
$D(2, 1 ; \alpha)$ into the Poisson superalgebra of pseudodifferential symbols on $S^{1|2}$.
We show that for the standard embedding of $D(2, 1 ; \alpha)$ into the Poisson superalgebra of differential operators on $S^{1|2}$ , the infinitesimal deformations correspond to formal deformations.
For the embedding of $D(2, 1 ; \alpha)$ into
the derived contact superconformal algebra ${K}'(4)$,
the infinitesimal deformations are formal deformations.}

\vskip 0.1in
\noindent
{\it JGP SC:} Lie superalgebras
\vskip 0.1in
\noindent

\vskip 0.5in
\noindent
{\bf 1. Introduction}

\vskip 0.5in

Recall that $D(2, 1 ; \alpha)$, where  $\alpha \in \C\backslash
\lbrace 0, -1 \rbrace$, is a one-parameter family of classical simple Lie superalgebras of dimension 17 [1].
The bosonic part of $D(2, 1 ; \alpha)$ is
$\s\l(2)\oplus \s\l(2)\oplus \s\l(2)$, and the action of
$D(2, 1 ; \alpha)_{\bar 0}$ on
$D(2, 1 ; \alpha)_{\bar 1}$
is the  product of two-dimensional representations.

An equivalent definition of this superalgebra was given in [2]. It is
the Lie superalgebra
$\Gamma (\sigma_1, \sigma_2,  \sigma_3)$, where $\sigma_i$ are nonzero complex numbers
such that $\sigma_1+\sigma_2+\sigma_3 = 0$.
Note that $\Gamma (\sigma_1, \sigma_2, \sigma_3)\cong D(2, 1 ; \alpha)$, where
$\alpha = \sigma_1/{\sigma_2}$.

In [3] we obtained
an embedding of  $D(2, 1 ; \alpha)$
into the Lie superalgebra of pseudodifferential symbols on $S^{1|2}$.
The integrability of infinitesimal
deformations of embeddings of Lie algebras
were  studied by Nijenhuis and Richardson in [4, 5].
The deformations of the standard embeddings
of $Vect(S^1)$
into the Poisson  algebra of
pseudodifferential symbols on $S^1$ and into the Lie algebra of
pseudodifferential symbols on $S^1$
were classified in [6, 7].
Similar problems for
the Lie superalgebras
$K(1)$ and $K(2)$
of contact vector fields
on the supercircles
$S^{1|1}$ and $S^{1|2}$ were studied in [8, 9].

In this work we consider the standard embedding of
$\Gamma_{\alpha} = \Gamma (2, -1 - \alpha,  \alpha - 1)$,
where $\alpha \in \C$,
into the Poisson superalgebra $P(4)$ of pseudodifferential symbols on the
supercircle $S^{1|2}$.
If $\alpha = 0$, we have
the natural
embedding of $\s\p\o (2|4)$ into $P(4)$.

According to the Richardson-Nijenhuis theory,
the infinitesimal deformations of this embedding are
classified by $H^1(\Gamma_{\alpha}, P(4))$.
We prove that this cohomology space is two dimensional and compute the corresponding cocycles.
$\Gamma_{\alpha}$ is naturally embedded into a
subsuperalgebra of $P(4)$, which is isomorphic to the contact superconformal algebra $K(4)$.
It consists of the functions
for which the corresponding Hamiltonian vector fields
commute with the Euler vector field.
Note that in this work, the realization of $K(4)$ inside $P(4)$ differs from the realization
that we obtained in [3], where $K(4)$ consists of the elements of $P(4)$
for which the corresponding Hamiltonian vector fields
commute with a {\it semi}-Euler vector field.

We  prove that the cohomology space
$H^1(\Gamma_{\alpha}, [K(4), K(4)])$
is one dimensional and that the infinitesimal deformations
are indeed the formal deformations of the embedding.

Note that $\Gamma_{\alpha}\subset P^+(4)$, where
$P^+(4)$ is a subsuperalgebra of $P(4)$ formed by  differential operators.
We prove that the cohomology space $H^1(\Gamma_{\alpha}, P^+(4))$ is one dimensional
and that the infinitesimal deformations
correspond to  formal deformations of the embedding.
The analogous results hold for the embedding of $\Gamma_{\alpha}$ into
the Lie superalgebra $P_{\hbox{h}}^+(4)$   of differential operators on $S^{1|2}$,
which contracts to $P^+(4)$.

\vskip 0.5in
\noindent
{\bf 2. Superalgebras $\Gamma(\sigma_1, \sigma_2, \sigma_3)$}

\vskip 0.5in

Recall the definition of $\Gamma(\sigma_1, \sigma_2, \sigma_3)$ [2].
Let $\g = \g_{\bar{0}} \oplus \g_{\bar{1}}$ be a Lie superalgebra, where
$\g_{\bar{0}} = sp(\psi_1)\oplus sp(\psi_2)\oplus sp(\psi_3)$ and
$\g_{\bar{1}} = V_1\otimes V_2\otimes V_3$, where
$V_i$ are two-dimensional vector spaces, and
$\psi_i$ is a non-degenerate skew-symmetric form on $V_i$, $i = 1, 2, 3$.
A representation of $\g_{\bar{0}}$ on $\g_{\bar{1}}$ is the tensor product
of the standard representations of $sp(\psi_i)$ in $V_i$.
Consider the $sp(\psi_i)$ - invariant bilinear mapping
$$\P_i: V_i\times V_i \rightarrow sp(\psi_i), \quad i = 1, 2, 3,$$
given by
$$\P_i(x_i, y_i)z_i = \psi_i(y_i, z_i)x_i - \psi_i(z_i, x_i)y_i$$
for all $x_i, y_i, z_i\in V_i$.
Let $\P$ be a mapping
$$\P:\g_{\bar{1}}\times \g_{\bar{1}}\rightarrow \g_{\bar{0}}$$
given by
\begin{equation*}
\begin{aligned}
&\P(x_1\otimes x_2\otimes x_3, y_1\otimes y_2\otimes y_3) =\\
&\sigma_1\psi_2(x_2, y_2)\psi_3(x_3, y_3)\P_1(x_1, y_1) + \\
&\sigma_2\psi_1(x_1, y_1)\psi_3(x_3, y_3)\P_2(x_2, y_2) + \\
&\sigma_3\psi_1(x_1, y_1)\psi_2(x_2, y_2)\P_3(x_3, y_3)
\end{aligned}
\end{equation*}
for all $x_i, y_i \in V_i, i = 1, 2, 3$,
where $\sigma_1, \sigma_2, \sigma_3$ are some complex numbers.
The super Jacobi identity is satisfied if and only if
$\sigma_1 + \sigma_2 + \sigma_3 = 0$. In this case $\g$ is denoted by
$\Gamma(\sigma_1, \sigma_2, \sigma_3)$.
Superalgebras $\Gamma(\sigma_1, \sigma_2, \sigma_3)$ and
$\Gamma(\sigma_1', \sigma_2', \sigma_3')$ are isomorphic if and only if there exists a nonzero element $k\in \C$ and a permutation $\pi$ of the set $\lbrace 1, 2, 3\rbrace$ such that
$$ \sigma_i' = k\cdot\sigma_{\pi i} \hbox{ for } i = 1, 2, 3.$$
Superalgebras $\Gamma(\sigma_1, \sigma_2, \sigma_3)$ are simple if
and only if $\sigma_1, \sigma_2, \sigma_3$ are all different from zero.
Note that $\Gamma(\sigma_1, \sigma_2, \sigma_3)\cong D(2, 1 ; \alpha)$ (see [1]) where
$\alpha = \sigma_1/{\sigma_2}$.

\vskip 0.5in
\noindent
{\bf 3. Embeddings  of $\Gamma(\sigma_1, \sigma_2, \sigma_3)$}

\vskip 0.5in

 The {\it Poisson algebra $P$ of pseudodifferential symbols on
the circle} is formed by the formal series
$$A(t, \tau) = \sum_{-\infty}^na_i(t)\tau^i,$$
where $a_i(t)\in \C[t, t^{-1}]$,
and the even variable $\tau$ corresponds to $\partial_t$, see [6].
The Poisson bracket is defined as follows:
$$\lbrace A(t, \tau), B(t, \tau)\rbrace = \partial_{\tau}A(t, \tau)\partial_tB(t, \tau) - \partial_tA(t, \tau)\partial_{\tau}B(t, \tau).$$

Let $\Lambda(2N)$ be the Grassmann algebra in $2N$ variables
$\xi_1, \ldots, \xi_N, \eta_1, \ldots, \eta_N $ with the parity $p(\xi_i) = p(\eta_i) = \bar{1}$.
The {\it Poisson superalgebra of pseudodifferential symbols on $S^{1|N}$} is
$P(2N) = P\otimes \Lambda(2N)$. The Poisson bracket is defined as follows:
$$\lbrace A, B\rbrace = \partial_{\tau}A\partial_tB - \partial_tA\partial_{\tau}B +
(-1)^{p(A)+1}\sum_{i = 1}^N(\partial_{\xi_i}A\partial_{\eta_i}B + \partial_{\eta_i}A\partial_{\xi_i}B).$$

Let
$P^+\subset P$ be the subalgebra of {\it differential operators}:
$$P^+ = \lbrace \sum_{i\geq 0}^na_i(t)\tau^i \rbrace.$$
Correspondingly,
$P^+(2N) = P^+\otimes \Lambda(2N)$.

Let
$W(2N)$ be the Lie superalgebra of all superderivations of the associative
superalgebra $\C [t, t^{-1}]\otimes \Lambda (2N)$.
By definition,
$$
K(2N) = \lbrace D \in W(2N)\mid D\Omega  = f\Omega \hbox{ for some }
f\in \C [t, t^{-1}]\otimes \Lambda (2N)\rbrace,
$$
where
$\Omega = dt + \sum_{i=1}^N \xi_id\eta_i + \eta_id\xi_i$
is a differential 1-form, which is called a {\it contact form} [10].
Note that there exists an embedding
$$K(2N)\subset P(2N),\quad N\geq 0.$$
Consider a $\Z$-grading $P(2N) = \oplus_i P_{(i)}P(2N)$
on the associative superalgebra $P(2N)$, defined by
\begin{equation*}
\begin{aligned}
&\hbox{deg } t = \hbox{deg } \eta_i =
\hbox{deg } \tau = \hbox{deg } \xi_i = 1 \hbox{ for } i = 1,
 \ldots, N.
 \end{aligned}
 \end{equation*}
With respect to the Poisson super-bracket,
$$\lbrace P_{(i)}(2N), P_{(j)}(2N)\rbrace \subset P_{({i+j-2})}(2N).$$
Thus $P_{(2)}(2N)$ is a subsuperalgebra of $P(2N)$, and
$P_{(2)}(2N) \cong K(2N)$.

\vskip 0.1in
\noindent
{\bf Remark 3.1.} To explain  this isomorphism, consider the analogous $\Z$-grading of the
associative algebra $P = \oplus_iP_{(i)}$ defined by the condition
$$\hbox{deg } t = \hbox{deg }\tau = 1.$$
Note that $P_{(2)}$ is a subalgebra of $P$,  which is isomorphic to the centerless Virasoro algebra $L = \oplus_nL_n$, where
$L_n = {1\over 2} t^{n+1}\tau^{-n +1}$, so that
$$[L_n, L_m] = (m - n)L_{n+m}.$$
Note that the elements of the Poisson algebra $P$ are functions $A(t, \tau)$
on the cylinder
$\dot{T}^*S^1 = T^*S^1\setminus S^1$, and they define the corresponding Hamiltonian vector fields on this manifold:
$$A(t, \tau) \longrightarrow H_A =
\partial_{\tau} A\partial_t - \partial_t A\partial_{\tau}.\eqno (3.1)$$
The subalgebra $P_{(2)}$ is formed precisely by such functions $A(t, \tau)$
for which  $H_A$
commutes with the Euler vector field:
$$[H_A, t\partial_t + \tau\partial_{\tau}] = 0.$$
The analogue of the formula (3.1) in the supercase is
as follows:
$$A(t, \tau,\xi_i, \eta_i) \longrightarrow H_A =
\partial_{\tau} A\partial_t - \partial_t A\partial_{\tau}
-(-1)^{p(A)}\sum_{i=1}^N(\partial_{\xi_i}A\partial_{\eta_i} +
\partial_{\eta_i} A\partial_{\xi_i}).$$
Then $K(2N)$ is defined as the set of all Hamiltonian functions
$A(t, \tau,\xi_i, \eta_i)$ for which  $H_A$
commutes with the Euler vector field:
$$[H_A, t\partial_t + \tau\partial_{\tau} +
\sum_{i=1}^N\xi_i\partial_{\xi_i} + \eta_i\partial_{\eta_i}] = 0.$$
Equivalently, $K(2N)$ is the subset of  $P(2N)$ of degree 2.

\noindent
Note that in [11] we considered a different
embedding of $K(2N)$ into $P(2N)$. It is defined by the condition that
$H_A$ commutes with a {\it semi}-Euler vector field on  $\dot{T}^*S^{1|N}$.

\vskip 0.1in
$K(2N)$ is simple if $N \not= 2$, and if $N = 2$, then the derived Lie superalgebra
$K'(4) = [K(4), K(4)]$
 is a simple ideal in $K(4)$
of codimension one, defined from the exact sequence
$$0\rightarrow K'(4)\rightarrow K(4)\rightarrow \C t^{-1}\tau^{-1}\xi_1\xi_2\eta_1\eta_2\rightarrow 0. $$

\vskip 0.1in
\noindent
{\bf Proposition 3.2.} For each $\alpha\in\C$
there exists an embedding
$$\rho_{\alpha}: \Gamma (2, -1 - \alpha, \alpha - 1)
\rightarrow K'(4)\subset P(4). $$
$\Gamma_{\alpha} = \rho_{\alpha}(\Gamma (2, -1 - \alpha, \alpha - 1))$ is
spanned by the following elements:
\begin{equation*}
\begin{aligned}
&E_{\alpha}^1 =  t^2, \quad
F_{\alpha}^1 =  \tau^2 - 2\alpha t^{-2}\xi_1\xi_2\eta_1\eta_2, \quad
H_{\alpha}^1 =  t\tau,\\
&E_{\alpha}^2 = \xi_1\xi_2,\quad
F_{\alpha}^2 = \eta_1\eta_2, \quad
H_{\alpha}^2 =  \xi_1\eta_1 + \xi_2\eta_2, \\
&E_{\alpha}^3 =  \xi_1\eta_2, \quad
F_{\alpha}^3 = \xi_2\eta_1,\quad
H_{\alpha}^3 =  \xi_1\eta_1 - \xi_2\eta_2,\\
&T_{\alpha}^1 =  t\eta_1, \quad
T_{\alpha}^2 =  t\eta_2 ,\quad
T_{\alpha}^3 =  t\xi_1, \quad
T_{\alpha}^4 =  t\xi_2,\\
&D_{\alpha}^1 =  \tau\xi_1 + \alpha t^{-1}\xi_1\xi_2\eta_2,\quad
D_{\alpha}^2 = \tau\xi_2 - \alpha t^{-1}\xi_1\xi_2\eta_1,\\
&D_{\alpha}^3 =  \tau\eta_1 + \alpha t^{-1}\xi_2\eta_1\eta_2,\quad
D_{\alpha}^4 =  \tau\eta_2 - \alpha t^{-1}\xi_1\eta_1\eta_2.
\end{aligned}
\tag{3.2}
\end{equation*}

\noindent
{\it Proof.}
Note that if $\alpha = 0$, then
$\Gamma (2, -1, - 1)\cong \s\p\o (2|4)$, and $\rho_{\alpha}$ is the standard
embedding of $\s\p\o (2|4)$ into $P(4)$.

Let
\begin{equation*}
\begin{aligned}
&V_1 =\hbox{Span}(e_1, e_2) , \quad
V_2   = \hbox{Span}(f_1, f_2), \quad
V_3 = \hbox{Span}(h_1, h_2),
\end{aligned}
\end{equation*}
and
\begin{equation*}
\begin{aligned}
&\psi_1(e_1, e_2) = - \psi_1 (e_2, e_1) = 1,\\
& \psi_2(f_1, f_2) = - \psi_2(f_2, f_1) = 1,\\
&\psi_3(h_1, h_2) = - \psi_3(h_2, h_1) = 1.
\end{aligned}
\end{equation*}
Explicitly an embedding $\rho_{\alpha}$ is given as follows:
\begin{equation*}
\begin{aligned}
&\rho_{\alpha}(\P_1(e_1, e_1)) = -E_{\alpha}^1, \quad
\rho_{\alpha}(\P_1(e_2, e_2)) = -F_{\alpha}^1,\quad
\rho_{\alpha}(\P_1(e_1, e_2)) = -H_{\alpha}^1,\\
&\rho_{\alpha}(\P_2(f_1, f_1)) = - 2F_{\alpha}^2, \quad
\rho_{\alpha}(\P_2(f_2, f_2)) = - 2E_{\alpha}^2,\quad
\rho_{\alpha}(\P_2(f_1, f_2)) =  H_{\alpha}^2,\\
&\rho_{\alpha}(\P_3(h_1, h_1)) = - 2F_{\alpha}^3, \quad
\rho_{\alpha}(\P_3(h_2, h_2)) =  2E_{\alpha}^3,\quad
\rho_{\alpha}(\P_3(h_1, h_2)) =  H_{\alpha}^3,\\
&\rho_{\alpha}(e_1\otimes f_1\otimes h_1) = {\sqrt 2}iT_{\alpha}^1, \quad
\rho_{\alpha}(e_1\otimes f_1\otimes h_2) = {\sqrt 2}iT_{\alpha}^2, \\
&\rho_{\alpha}(e_1\otimes f_2\otimes h_1) = -{\sqrt 2}iT_{\alpha}^4, \quad
\rho_{\alpha}(e_1\otimes f_2\otimes h_2) = {\sqrt 2}iT_{\alpha}^3,\\
&\rho_{\alpha}(e_2\otimes f_1\otimes h_1) = {\sqrt 2}iD_{\alpha}^3, \quad
\rho_{\alpha}(e_2\otimes f_1\otimes h_2) = {\sqrt 2}iD_{\alpha}^4, \\
&\rho_{\alpha}(e_2\otimes f_2\otimes h_1) = -{\sqrt 2}iD_{\alpha}^2, \quad
\rho_{\alpha}(e_2\otimes f_2\otimes h_2) = {\sqrt 2}iD_{\alpha}^1. \\
\end{aligned}
\end{equation*}
Thus $sp(\psi_i) \cong \hbox{Span} (E_{\alpha}^i, H_{\alpha}^i, F_{\alpha}^i)$
for $i = 1, 2, 3$.
$$\eqno\Q$$

\vskip 0.5in
\noindent
{\bf 4. Deformations of embeddings}

\vskip 0.5in

Let $\rho:\g\rightarrow \h$ be an embedding of Lie superalgebras,
then $\h$ is a $\g$-module. A map $\rho + \beta\rho_1: \g\rightarrow \h$, where
$\rho_1\in Z^1(\g, \h)$ is a Lie superalgebra homomorphism up to quadratic terms in $\beta$.
It is called an infinitesimal deformation.
Infinitesimal deformations are classified by $H^1(\g, \h)$, see [4, 5].

Let
$$\tilde{\rho}_{\beta} = \rho + \sum_{k = 1}^{\infty}{\beta}^k\rho_k:\g\rightarrow \h,$$
where $\rho_k:\g\rightarrow \h$ are even linear maps,
satisfy
$$\tilde{\rho}_{\beta}([X, Y]) = [\tilde{\rho}_{\beta}(X), \tilde{\rho}_{\beta}(Y)].$$
$\tilde{\rho}_{\beta}$ is called a formal deformation of $\rho$.
The integrability conditions are conditions for existence of formal deformation corresponding to a given infinitesimal deformation. The obstructions for existence of a formal deformation belong to
the second cohomology group $H^2(\g, \h)$, see [4, 6].

Let $\varphi_{\beta} = \tilde{\rho}_{\beta} - \rho$. Then
$$[\varphi_{\beta}(X), \rho(Y)] +  [\rho(X), \varphi_{\beta}(Y)] -  \varphi_{\beta}([X, Y]) +
\sum_{i, j> 0}[\rho_i(X),\rho_j(Y)]{\beta}^{i+j} = 0. \eqno (4.1)$$

The first three terms are $(d\varphi_{\beta})(X, Y)$, where $d$ stands for coboundary.
For arbitrary linear maps $\varphi, \varphi': \g\rightarrow\h$, define

\begin{equation*}
\begin{aligned}
&\lbrack \lbrack\varphi, \varphi'\rbrack \rbrack: \g\otimes\g\rightarrow \h,\\
&\lbrack \lbrack\varphi, \varphi'\rbrack \rbrack (X, Y) =
\lbrack\varphi(X), \varphi'(Y)\rbrack + \lbrack\varphi'(X), \varphi(Y)\rbrack.\\
\end{aligned}
\tag{4.2}
\end{equation*}
Relation (4.1) is equivalent to
$$d\varphi_{\beta} + {1\over 2}\lbrack \lbrack\varphi_{\beta}, \varphi_{\beta}\rbrack \rbrack = 0.$$
Expanding this relation in power series in $\beta$, we have
$$d\rho_k + {1\over 2}\sum_{i + j = k}\lbrack\lbrack\rho_i,\rho_j\rbrack\rbrack = 0.$$
The first nontrivial relation is
$$d\rho_2 + {1\over 2}\lbrack\lbrack\rho_1,\rho_1\rbrack\rbrack = 0,$$
and it gives the first obstruction to integrability of an infinitesimal deformation.
Note that (4.2) defines  a bilinear map, called the cup product:
$$H^1(\g, \h)\otimes H^1(\g, \h)\rightarrow H^2(\g, \h).$$
The obstructions to integrability of infinitesimal deformations lie in $H^2(\g, \h)$.
Thus we have to compute $H^1(\g, \h)$ and the product classes in $H^2(\g, \h)$.

\noindent
Consider the embedding (3.2).

\vskip 0.1in
\noindent
{\bf Theorem 4.1.}
$\hbox{dim}H^1(\Gamma_{\alpha}, P(4)) = 2$. The cohomology space is spanned by the classes of the 1-cocycles $\theta_1$
and $\theta_2$ given as follows:

\begin{equation*}
\begin{aligned}
&\theta_1(D_{\alpha}^1) =  t^{-1}\xi_1,\quad
\theta_1(D_{\alpha}^2) =&  t^{-1}\xi_2,\\
&\theta_1(D_{\alpha}^3) =  t^{-1}\eta_1,\quad
\theta_1(D_{\alpha}^4) =&  t^{-1}\eta_2,\\
&\theta_1(F_{\alpha}^1) = 2t^{-1}\tau,\quad
\theta_1(H_{\alpha}^1) =& 1.\\
\end{aligned}
\tag{4.3}
\end{equation*}

\begin{equation*}
\begin{aligned}
&\theta_2(T_{\alpha}^3) =
{\tau}^{-1}\xi_1-t^{-1}{\tau}^{-2}
\xi_1\xi_2\eta_2,\quad
\theta_2(T_{\alpha}^4) =
{\tau}^{-1}\xi_2+t^{-1}{\tau}^{-2}
\xi_1\xi_2\eta_1,\\
&\theta_2(D_{\alpha}^1) =  t^{-1}\xi_1,\quad
\theta_2(D_{\alpha}^2) =  t^{-1}\xi_2,\\
&\theta_2(D_{\alpha}^3) =
-(1+\alpha) t^{-2}{\tau}^{-1}\xi_2\eta_1\eta_2,\quad
\theta_2(D_{\alpha}^4) =
(1+\alpha) t^{-2}{\tau}^{-1}\xi_1\eta_1\eta_2,\\
&\theta_2(E_{\alpha}^1) =
t\tau^{-1}-\tau^{-2}\xi_1\eta_1
-\tau^{-2}\xi_2\eta_2
-2t^{-1}\tau^{-3}\xi_1\xi_2\eta_1\eta_2,\quad
\theta_2(E_{\alpha}^2)=
t^{-1}\tau^{-1}\xi_1\xi_2,\\
&\theta_2(F_{\alpha}^1) =
t^{-1}\tau + t^{-2}\xi_1\eta_1
+ t^{-2}\xi_2\eta_2
+2(1+\alpha)t^{-3}\tau^{-1}
\xi_1\xi_2\eta_1\eta_2, \\
&\theta_2(F_{\alpha}^2)=
-t^{-1}\tau^{-1}\eta_1\eta_2,\quad
\theta_2(H_{\alpha}^1) = 1.\\
\end{aligned}
\end{equation*}

\noindent
{\it Proof.}
Consider $\g\l(2) \cong  \hbox{Span}(\xi_i\eta_j \hbox{ }|\hbox{ } i, j = 1, 2)\subset \Gamma_{\alpha}$.
The diagonal subalgebra of $\g\l(2)$ consists of
$h = h_1\xi_1\eta_1 + h_2\xi_2\eta_2$, where $h_1, h_2\in \C$.
Let $\epsilon_i(h) = h_i, i = 1, 2$. Obviously, $\hbox{Span}(\xi_1, \xi_2)$ is the standard $\g\l(2)$-module,
$\hbox{Span}(\eta_1, \eta_2)$ is its dual,
$\xi_i$ and $\eta_i$ have weights $\epsilon_i$ and $-\epsilon_i$. Note that
$H^1(\Gamma_{\alpha}, P(4))$ is a trivial $\g\l(2)$-module, since a Lie (super)algebra
acts trivially on its cohomology [12]. Hence we have to compute only the $1$-cocycles of weight zero.
Also note that
$$H^1(\Gamma_{\alpha}, P(4)) = \oplus_{k\in\Z}H^{1}(\Gamma_{\alpha}, P_{(k)}(4)),
\eqno (4.4)$$
and
$$H^1(\Gamma_{\alpha}, P(4)_{(k)}) = \oplus_{n\in\Z}H^{1,n}(\Gamma_{\alpha}, P_{(k)}(4)),
\eqno (4.5)$$
where the $\Z$-grading is given by the condition
$$\hbox{deg } t = 1, \hbox{ deg } \tau = -1,\quad \hbox{deg }\xi_i  =  \hbox{deg }\eta_i = 0.$$

\noindent
Let $c\in C^{1,n}(\Gamma_{\alpha}, P_{(k)}(4))$
be a $1$-cochain of weight zero.

Assume that $k$ is even: $k = 2k_1$.
If $c\not= 0$, then $n$
must be even: $n =2m$, and
 $c$ acts on the odd elements of $\Gamma_{\alpha}$ as follows:

\begin{equation*}
\begin{aligned}
&c(T_{\alpha}^1) = g_1^mt^{m+k_1}\tau^{-m+k_1-1}\eta_1
+ s_1^mt^{m+k_1-1}\tau^{-m+k_1-2}\xi_2\eta_1\eta_2,\\
&c(T_{\alpha}^2) = g_2^mt^{m+k_1}\tau^{-m+k_1-1}\eta_2 +
s_2^mt^{m+k_1-1}\tau^{-m+k_1-2}\xi_1\eta_1\eta_2,\\
&c(T_{\alpha}^3) = g_3^{m}t^{m+k_1}\tau^{-m+k_1-1}\xi_1 +
s_3^mt^{m+k_1-1}\tau^{-m+k_1-2}\xi_1\xi_2\eta_2, \\
&c(T_{\alpha}^4) = g_4^mt^{m+k_1}\tau^{-m+k_1-1}\xi_2 +
s_4^mt^{m+k_1-1}\tau^{-m+k_1-2}\xi_1\xi_2\eta_1,\\
&c(D_{\alpha}^1) = r_1^mt^{m+k_1-1}\tau^{-m+k_1}\xi_1 +
q_1^mt^{m+k_1-2}\tau^{-m+k_1-1}\xi_1\xi_2\eta_2,\\
&c(D_{\alpha}^2) = r_2^mt^{m+k_1-1}\tau^{-m+k_1}\xi_2 +
q_2^mt^{m+k_1-2}\tau^{-m+k_1-1}\xi_1\xi_2\eta_1,\\
&c(D_{\alpha}^3) = r_3^mt^{m+k_1-1}\tau^{-m+k_1}\eta_1 +
q_3^mt^{m+k_1-2}\tau^{-m+k_1-1}\xi_2\eta_1\eta_2,\\
&c(D_{\alpha}^4) = r_4^mt^{m+k_1-1}\tau^{-m+k_1}\eta_2 +
q_4^mt^{m+k_1-2}\tau^{-m+k_1-1}\xi_1\eta_1\eta_2,\\
\end{aligned}
\tag{4.6}
\end{equation*}
where $g_i^m, s_i^m, r_i^m, q_i^m\in\C$.
Let
\begin{equation*}
\begin{aligned}
&c_0 = t^{m+k_1}\tau^{-m+k_1}, \quad c_1 = t^{m+k_1-1}\tau^{-m+k_1-1}\xi_1\eta_1,\\
&c_2 = t^{m+k_1-1}\tau^{-m+k_1-1}\xi_2\eta_2, \quad
c_3 = t^{m+k_1-2}\tau^{-m+k_1-2}\xi_1\xi_2\eta_1\eta_2.\\
\end{aligned}
\end{equation*}
The elements of weight zero in
$C^{0,2m}(\Gamma_{\alpha}, P_{(k)}(4))$ span the subspace $\hbox{Span}(c_0, c_1, c_2, c_3)$.
Note that the coefficients $g_i^m$ in (4.6) are as follows:

\begin{equation*}
\begin{aligned}
&\hbox{if } c = dc_0, \hbox{ then } g_1^m = g_2^m = g_3^m = g_4^m = m - k_1,\\
&\hbox{if } c = dc_1, \hbox{ then } g_1^m = -g_3^m = 1, g_2^m = g_4^m = 0,\\
&\hbox{if } c = dc_2, \hbox{ then } g_1^m = g_3^m = 0, g_2^m = -g_4^m = 1.\\
\end{aligned}
\tag{4.7}
\end{equation*}

Assume that $k$ is odd: $k = 2k_1+1$.
If $c\not= 0$, then $n$
must be odd: $n =2m+1$, and
 $c$ acts on the odd elements of $\Gamma_{\alpha}$ as follows:

\begin{equation*}
\begin{aligned}
&c(T_{\alpha}^1) = g_1^mt^{m+k_1+1}\tau^{-m+k_1-1}\eta_1
+ s_1^mt^{m+k_1}\tau^{-m+k_1-2}\xi_2\eta_1\eta_2,\\
&c(T_{\alpha}^2) = g_2^mt^{m+k_1+1}\tau^{-m+k_1-1}\eta_2 +
s_2^mt^{m+k_1}\tau^{-m+k_1-2}\xi_1\eta_1\eta_2,\\
&c(T_{\alpha}^3) = g_3^{m}t^{m+k_1+1}\tau^{-m+k_1-1}\xi_1 +
s_3^mt^{m+k_1}\tau^{-m+k_1-2}\xi_1\xi_2\eta_2, \\
&c(T_{\alpha}^4) = g_4^mt^{m+k_1+1}\tau^{-m+k_1-1}\xi_2 +
s_4^mt^{m+k_1}\tau^{-m+k_1-2}\xi_1\xi_2\eta_1,\\
&c(D_{\alpha}^1) = r_1^mt^{m+k_1}\tau^{-m+k_1}\xi_1 +
q_1^mt^{m+k_1-1}\tau^{-m+k_1-1}\xi_1\xi_2\eta_2,\\
&c(D_{\alpha}^2) = r_2^mt^{m+k_1}\tau^{-m+k_1}\xi_2 +
q_2^mt^{m+k_1-1}\tau^{-m+k_1-1}\xi_1\xi_2\eta_1,\\
&c(D_{\alpha}^3) = r_3^mt^{m+k_1}\tau^{-m+k_1}\eta_1 +
q_3^mt^{m+k_1-1}\tau^{-m+k_1-1}\xi_2\eta_1\eta_2,\\
&c(D_{\alpha}^4) = r_4^mt^{m+k_1}\tau^{-m+k_1}\eta_2 +
q_4^mt^{m+k_1-1}\tau^{-m+k_1-1}\xi_1\eta_1\eta_2,\\
\end{aligned}
\tag{4.8}
\end{equation*}
where $g_i^m, s_i^m, r_i^m, q_i^m\in\C$.

The elements of weight zero in
$C^{0,2m+1}(\Gamma_{\alpha}, P_{(k)}(4))$
span the subspace $\hbox{Span}(c_0, c_1, c_2, c_3)$,
where
\begin{equation*}
\begin{aligned}
&c_0 = t^{m+k_1+1}\tau^{-m+k_1}, \quad c_1 = t^{m+k_1}\tau^{-m+k_1-1}\xi_1\eta_1,\\
&c_2 = t^{m+k_1}\tau^{-m+k_1-1}\xi_2\eta_2, \quad
c_3 = t^{m+k_1-1}\tau^{-m+k_1-2}\xi_1\xi_2\eta_1\eta_2,\\
\end{aligned}
\tag{4.9}
\end{equation*}
and the coefficients $g_i^m$ are as in (4.7).

Let $X, Y\in \Gamma_{\alpha}$.
Note that
\begin{equation*}
\begin{aligned}
&dc(X, Y) = \lbrace X, c(Y)\rbrace +  \lbrace Y, c(X)\rbrace - c(\lbrace X, Y\rbrace),
\quad \hbox{ if } p(X) = p(Y) = \bar{1},\\
&dc(X, Y) = \lbrace X, c(Y)\rbrace -  \lbrace Y, c(X)\rbrace - c(\lbrace X, Y\rbrace),
\quad \hbox{ if } p(X) = \bar{0}, \hbox{  }  p(Y) = \bar{1},\\
&dc(X, Y) = \lbrace X, c(Y)\rbrace -  \lbrace Y, c(X)\rbrace - c(\lbrace X, Y\rbrace),
\quad \hbox{ if } p(X) =  p(Y) = \bar{0}.\\
\end{aligned}
\end{equation*}
Let $c\in Z^{1,n}(\Gamma_{\alpha}, P_{(k)}(4))$ be of weight zero.
The condition $dc(X, Y) = 0$ gives that
\begin{equation*}
\begin{aligned}
&\lbrace T_{\alpha}^1, c(T_{\alpha}^3)\rbrace +  \lbrace T_{\alpha}^3, c(T_{\alpha}^1)\rbrace - c(E_{\alpha}^1) = 0,\\
&\lbrace T_{\alpha}^2, c(T_{\alpha}^4)\rbrace +
\lbrace T_{\alpha}^4, c(T_{\alpha}^2)\rbrace - c(E_{\alpha}^1) = 0.
\end{aligned}
\tag{4.10}
\end{equation*}
It follows that
$$g_3^m + g_1^m = g_2^m + g_4^m. \eqno (4.11)$$

\noindent
{\it Case }$m\not= k_1$.
One can change $c$ by adding (or removing) coboundaries $dc_i$ for $i = 0, 1, 2$,
and thus
assume that $g_i^m = 0$ for
$i = 1, 2, 3$. Then from (4.11) $g_4^m = 0$.
 Note that
$$\lbrace T_{\alpha}^1, c(T_{\alpha}^2)\rbrace +
\lbrace T_{\alpha}^2, c(T_{\alpha}^1)\rbrace = 0,\eqno (4.12)$$
hence $s_2^m = -s_1^m$.
$$\lbrace T_{\alpha}^1, c(T_{\alpha}^4)\rbrace +
\lbrace T_{\alpha}^4, c(T_{\alpha}^1)\rbrace = 0, \eqno (4.13)$$
hence $s_4^m = -s_1^m$.
$$\lbrace T_{\alpha}^2, c(T_{\alpha}^3)\rbrace +
\lbrace T_{\alpha}^3, c(T_{\alpha}^2)\rbrace = 0, \eqno (4.14)$$
hence $s_3^m = -s_2^m = s_1^m.$

Note that if $c = dc_3$, then in (4.6) $g_i^m = 0$ for $i = 1, \dots, 4$ and
$s_1^m = s_3^m = -s_2^m = -s_4^m = 1$.
Changing $1$-cocycle $c$ using the coboundary $dc_3$, we can assume
in addition to all $g_i^m = 0$ that $s_i^m = 0$
for $i = 1, 2, 3,  4$.
Hence $c(T_{\alpha}^i) = 0$ for
$i = 1, 2, 3,  4$.
We have that
$$\lbrace E_{\alpha}^2, c(T_{\alpha}^1)
\rbrace -  \lbrace T_{\alpha}^1, c(E_{\alpha}^2)\rbrace + c(T_{\alpha}^4)
= 0,\eqno(4.15)$$
hence $c(E_{\alpha}^2) = 0$.
$$\lbrace E_{\alpha}^3, c(T_{\alpha}^1)
\rbrace -  \lbrace T_{\alpha}^1, c(E_{\alpha}^3)
\rbrace -  c(T_{\alpha}^2) = 0, \eqno(4.16)$$
hence $c(E_{\alpha}^3) = 0.$
Also
$$\lbrace F_{\alpha}^2, c(T_{\alpha}^3)
\rbrace -  \lbrace
T_{\alpha}^3, c(F_{\alpha}^2)\rbrace + c(T_{\alpha}^2) = 0, \eqno (4.17)$$
hence $c(F_{\alpha}^2) = 0$.
$$\lbrace F_{\alpha}^3, c(T_{\alpha}^3)
\rbrace -  \lbrace T_{\alpha}^3, c(F_{\alpha}^3)
\rbrace -  c(T_{\alpha}^4) = 0, \eqno (4.18)$$
hence $c(F_{\alpha}^3) = 0.$
Then
$$\lbrace D_{\alpha}^1, c(T_{\alpha}^4)\rbrace +
\lbrace T_{\alpha}^4, c(D_{\alpha}^1)\rbrace  = 0,\eqno  (4.19)$$
$$\lbrace T_{\alpha}^2, c(D_{\alpha}^1)
\rbrace +  \lbrace D_{\alpha}^1, c(T_{\alpha}^2)\rbrace  = 0. \eqno (4.20)$$
From (4.19)  $(-m+k_1)r_1^m + q_1^m = 0$, and from (4.20)
$(-m+k_1)r_1^m - q_1^m = 0$.
Hence,  $r_1^m = q_1^m = 0$.
$$\lbrace T_{\alpha}^1, c(D_{\alpha}^2)\rbrace +
\lbrace D_{\alpha}^2, c(T_{\alpha}^1)\rbrace  = 0, \eqno (4.21)$$
$$\lbrace D_{\alpha}^2, c(T_{\alpha}^3)\rbrace +
 \lbrace T_{\alpha}^3, c(D_{\alpha}^2)\rbrace  = 0. \eqno (4.22)$$
From (4.21) $(-m+k_1)r_2^m + q_2^m = 0$, and from (4.22) $(-m+k_1)r_2^m - q_2^m = 0$.
Hence,  $r_2^m = q_2^m = 0$.
$$\lbrace T_{\alpha}^2, c(D_{\alpha}^3)
\rbrace +  \lbrace D_{\alpha}^3, c(T_{\alpha}^2)\rbrace  = 0. \eqno (4.23)$$
$$\lbrace T_{\alpha}^4, c(D_{\alpha}^3)\rbrace +
\lbrace D_{\alpha}^3, c(T_{\alpha}^4)\rbrace  = 0.\eqno (4.24)$$
From (4.23) $(-m+k_1)r_3^m + q_3^m = 0$ and from (4.24)
$(-m+k_1)r_3^m - q_3^m = 0$.
Hence,  $r_3^m = q_3^m = 0$.
$$\lbrace T_{\alpha}^3, c(D_{\alpha}^4)\rbrace +
\lbrace D_{\alpha}^4, c(T_{\alpha}^3)\rbrace  = 0. \eqno (4.25)$$
$$\lbrace T_{\alpha}^1, c(D_{\alpha}^4)\rbrace +
\lbrace D_{\alpha}^4, c(T_{\alpha}^1)\rbrace  = 0,\eqno (4.26)$$
from (4.25) $(-m+k_1)r_4^m + q_4^m = 0$ and from (4.26) $(-m+k_1)r_4^m - q_4^m = 0$.
Hence,  $r_4^m = q_4^m = 0$.
Thus $c(D_{\alpha}^i) = 0$ for
$i = 1, 2, 3, 4$.

\noindent
Note that if $\alpha \not= \pm 1$, then $\Gamma_{\alpha}$ is simple. Hence
${(\Gamma_{\alpha})}_{\bar 0} =
[{(\Gamma_{\alpha})}_{\bar 1}, {(\Gamma_{\alpha})}_{\bar 1}]$.
Since $c({(\Gamma_{\alpha})}_{\bar 1})  = 0$, $c({(\Gamma_{\alpha})}_{\bar 0})  = 0$.
If $\alpha = \pm 1$, then
$$[{(\Gamma_{\alpha})}_{\bar 1}, {(\Gamma_{\alpha})}_{\bar 1}] =
\hbox{Span}(E_{\alpha}^i, H_{\alpha}^i, F_{\alpha}^i), \hbox{ where } i = 1, 2
\hbox{ or } i = 1, 3.$$
On the other hand, we proved that
$c(E_{\alpha}^i) = c(F_{\alpha}^i) = 0 \hbox { for } i = 2, 3$. Then it follows from
$\lbrace F_{\alpha}^i, E_{\alpha}^i \rbrace = H_{\alpha}^i$ that
$c(H_{\alpha}^i) = 0$ for  $i = 2, 3$. Hence $c({(\Gamma_{\alpha})}_{\bar 0})  = 0$, and
the cocycle $c$ vanishes.

\vskip 0.1in
\noindent
{\it Case }$m = k_1$.
Replacing  $1$-cocycle $c$ by  $c - d(c_1  + c_2)$,
we can assume that $g_1^m = g_2^m = 0$.
Then from (4.11) $g_3^m = g_4^m$.
In addition, changing $c$ by a multiple of $dc_3$,
 we can assume that $s_1^m = 0$.
From (4.12) $s_2^m = -s_1^m = 0$.
Hence $c(T_{\alpha}^1) = c(T_{\alpha}^2) = 0$.
Next, from
$$\lbrace T_{\alpha}^3, c(T_{\alpha}^4)\rbrace +
\lbrace T_{\alpha}^4, c(T_{\alpha}^3)\rbrace  = 0,$$
$s_3^m = -s_4^m$, and
from (4.13) $s_4^m = g_4^m$.
It follows from (4.10) that if $k$ is even ($k = 2m$), then
$$c(E_{\alpha}^1) = g_4^m(t^{2m+1}\tau^{-1} - t^{2m}\tau^{-2}\xi_1\eta_1
- t^{2m}\tau^{-2}\xi_2\eta_2 - 2t^{2m-1}\tau^{-3}\xi_1\xi_2\eta_1\eta_2),$$
and if $k$ is odd ($k = 2m+1$), then
$$c(E_{\alpha}^1) = g_4^m(t^{2m+2}\tau^{-1} - t^{2m+1}\tau^{-2}\xi_1\eta_1
- t^{2m+1}\tau^{-2}\xi_2\eta_2 - 2t^{2m}\tau^{-3}\xi_1\xi_2\eta_1\eta_2).$$
It follows from
$$\lbrace E_{\alpha}^1, c(D_{\alpha}^1)\rbrace -
\lbrace D_{\alpha}^1, c(E_{\alpha}^1)\rbrace  + 2c(T_{\alpha}^3) = 0 \eqno(4.27)$$
that
\begin{equation*}
\begin{aligned}
&2mg_4^mt^{2m}\tau^{-1}\xi_1 = 0, \hbox{ if } k \hbox{ is even, }\\
&(2m+1)g_4^mt^{2m+1}\tau^{-1}\xi_1 = 0, \hbox{ if } k \hbox{ is odd. }\\
\end{aligned}
\end{equation*}
Hence $g_4^m = 0$
(unless $m = k = 0$; we will consider this case separately).
Thus $c(T_{\alpha}^i) = 0$ for $i = 1, 2, 3, 4$, and $c(E_{\alpha}^1) = 0$.
From (4.15)
$c(E_{\alpha}^2) = 0$. From (4.16)
$c(E_{\alpha}^3) = 0$. From (4.17)
$c(F_{\alpha}^2) = 0$. From (4.18)
$c(F_{\alpha}^3) = 0$.
From (4.19)
$q_1^m = 0$.
From (4.21)
$q_2^m = 0$. From (4.23)
$q_3^m = 0$. From (4.25)
$q_4^m = 0$.

\noindent
From
$$\lbrace F_{\alpha}^3, c(D_{\alpha}^1)\rbrace -
\lbrace D_{\alpha}^1, c(F_{\alpha}^3)\rbrace  - c(D_{\alpha}^2) = 0,\eqno (4.28)$$
$r_1^m = r_2^m$. From
$$\lbrace F_{\alpha}^3, c(D_{\alpha}^4)\rbrace -
\lbrace D_{\alpha}^4, c(F_{\alpha}^3)\rbrace  + c(D_{\alpha}^3) = 0,\eqno (4.29)$$
$r_3^m = r_4^m$. From
$$\lbrace F_{\alpha}^2, c(D_{\alpha}^1)\rbrace -
\lbrace D_{\alpha}^1, c(F_{\alpha}^2)\rbrace  + c(D_{\alpha}^4) = 0,\eqno (4.30)$$
$r_1^m = r_4^m$. Thus
$r_1^m = r_2^m = r_3^m = r_4^m$.
One can easily see that $c$ is a multiple of $dc_0$.

\noindent
{\it Case: $m = k = 0$.}
We have already proved that
$c(T_{\alpha}^1)= c(T_{\alpha}^2) = 0$ and
\begin{equation*}
\begin{aligned}
&c(T_{\alpha}^3)=
g_4^0\tau^{-1}\xi_1 - g_4^0t^{-1}\tau^{-2}\xi_1\xi_2\eta_2,\\
&c(T_{\alpha}^4)=
g_4^0\tau^{-1}\xi_2 + g_4^0t^{-1}\tau^{-2}\xi_1\xi_2\eta_1.\\
\end{aligned}
\end{equation*}
Set $g_4^0 = 1$.
From (4.10)
$$c(E_{\alpha}^1) =
t\tau^{-1}-\tau^{-2}\xi_1\eta_1
-\tau^{-2}\xi_2\eta_2
-2t^{-1}\tau^{-3}\xi_1\xi_2\eta_1\eta_2,$$
from (4.16)
$c(E_{\alpha}^3)= 0$.
From (4.17)
$c(F_{\alpha}^2)= -t^{-1}\tau^{-1}\eta_1\eta_2$,
from (4.18)
$c(F_{\alpha}^3)= 0$.

\noindent
From (4.27) $q_1^0 = 0$,
from (4.21) $q_2^0 = 0$,
from
$$\lbrace T_{\alpha}^2, c(D_{\alpha}^3)\rbrace +
\lbrace D_{\alpha}^3, c(T_{\alpha}^2)\rbrace  -(1+\alpha)c(F_{\alpha}^2) = 0,$$
$q_3^0 = -1 - \alpha$.
From (4.25) $q_4^0 = \alpha + 1$.
From (4.28) $r_1^0 = r_2^0$,
from (4.29) $r_3^0 = r_4^0$,
and from (4.30) $r_1^0 = r_4^0 + 1$.

\noindent
Set $r_4^0 = 0$. Then we obtain the cocycle $\theta_2$.
Finally, $c - \theta_2 = r_4^0\theta_1$.
$$\eqno\Q$$

\vskip 0.1in
\noindent
{\bf Corollary 4.2.}
$\hbox{dim}H^1(\Gamma_{\alpha}, K'(4)) = 1$. The cohomology space is spanned by the class of
the $1$-cocycle $\theta$ given as follows:

\begin{equation*}
\begin{aligned}
&\theta(T_{\alpha}^1) = \tau^{-1}\xi_2\eta_1\eta_2,\quad\quad\quad\qquad
\theta(T_{\alpha}^2) = -\tau^{-1}\xi_1\eta_1\eta_2,\\
&\theta(T_{\alpha}^3) = \tau^{-1}\xi_1\xi_2\eta_2,\quad\quad\quad\qquad
\theta(T_{\alpha}^4) = -\tau^{-1}\xi_1\xi_2\eta_1,\\
&\theta(D_{\alpha}^1) =  t^{-1}\xi_1\xi_2\eta_2,\quad\qquad\qquad
\theta(D_{\alpha}^2) =  -t^{-1}\xi_1\xi_2\eta_1,\\
&\theta(D_{\alpha}^3) =  t^{-1}\xi_2\eta_1\eta_2,\qquad\hbox{ }\hbox{ }\qquad
\theta(D_{\alpha}^4) =  -t^{-1}\xi_1\eta_1\eta_2,\\
&\theta(E_{\alpha}^1) = 2\tau^{-2}\xi_1\xi_2\eta_1\eta_2,\hbox{ }\qquad\quad
\theta(F_{\alpha}^1) = -2t^{-2}\xi_1\xi_2\eta_1\eta_2.\\
\end{aligned}
\tag{4.31}
\end{equation*}
The map $\tilde{\rho}_{\alpha, \beta} = \rho_{\alpha} + \beta\theta$
($\beta\in\C$)  is a formal deformation of the embedding $\rho_{\alpha}$.

\noindent
{\it Proof}.
We already noticed that
$$P_{(2)}(4) \cong K(4).$$
From Theorem 4.1, $H^1(\Gamma_{\alpha}, K(4)) = 0$.
Hence $H^1(\Gamma_{\alpha}, K'(4)) = H^{1,0}(\Gamma_{\alpha}, K'(4))$ and it is spanned by
the class of the 1-cocycle
$dc_3$ for $m = 0$ and $k_1 = 1$, because
$c_3 = t^{-1}\tau^{-1}\xi_1\xi_2\eta_1\eta_2\not\in K'(4)$.
The coefficients in (4.6) for this cocycle   are $g_i^0 = r_i^0 = 0$, $s_1^0 = s_3^0 = -s_2^0 = -s_4^0 = 1$, and $q_1^0 = q_3^0 = -q_2^0 = -q_4^0 = 1$.
Thus $dc_3 = \theta$ as  is given in (4.31).

According to the Richardson-Nijenhuis theory, one has to determine the cup product
$\lbrack\lbrack \theta, \theta \rbrack\rbrack$ [6].
It is easy to see  that this cup product is identically zero (and not only
in cohomology). Thus $\tilde{\rho}_{\alpha, \beta} = \rho_{\alpha} +
\beta\theta$ is a formal deformation of the embedding $\rho_{\alpha}$.
$$\eqno\Q$$

Note that under the embedding (3.2) $\Gamma (2, -1-{\alpha}, {\alpha}-1)$ is
realized in differential operators. Thus
$$\rho_{\alpha}: \Gamma (2, -1-{\alpha}, {\alpha}-1) \longrightarrow P^+(4).\eqno (4.32)$$

\vskip 0.1in
\noindent
{\bf Theorem 4.3.}
$\hbox{dim}H^1(\Gamma_{\alpha}, P^+(4)) = 1$. The cohomology space is spanned by the class of the 1-cocycle
$\theta_1$ given in (4.3). The infinitesimal deformation defined by $\theta_1$
corresponds to a formal deformation of the embedding (4.32).

\noindent
{\it Proof}.
Let $c\in Z^{1,n}(\Gamma_{\alpha}, P_{(k)}^+(4))$
be a $1$-cocycle of weight zero.
Recall that
$k = 2k_1$ and
$n = 2m$ or
$k = 2k_1+1$ and
$n = 2m+1$.
It follows from (4.6) and (4.8)
that if $k_1 - m \geq 2$
or $k_1 - m < 0$, then
$$H^{1,n}(\Gamma_{\alpha}, P_{(k)}^+(4)) =
H^{1,n}(\Gamma_{\alpha}, P_{(k)}(4)) = 0.$$
If $k_1 - m = 1$, then $s_i^m = 0$
for $i = 1, 2, 3, 4$.
As in the proof of the case
$m\not= k_1$ of Theorem 4.1,
we have that $c$ has the zero cohomology class.

Finally, if $k_1 - m = 0$, then
in (4.6) and (4.8)
$g_i^m = s_i^m = q_i^m = 0$
for $i = 1, 2, 3, 4$.
As in the proof of the case
$m = k_1$ of Theorem 4.1, we obtain that $r_1^m =  r_2^m = r_3^m = r_4^m$.
Then if $k$ is odd or $k\not = 0$ is even, then $c$ is a multiple of $dc_0$.
If $k = m = 0$, then $dc_0 = 0$,
and $c = r_4^0\theta_1$.

\noindent
Let $\rho_1 = \theta_1$.
We will find a formal deformation corresponding to the infinitesimal deformation $\rho_{\alpha} + \beta\rho_1$.
 Recall that
$$d\rho_2  = - {1\over 2}\lbrack \lbrack\rho_1, \rho_1\rbrack \rbrack,\eqno (4.33)$$
where $\rho_2$ is a linear map
$$\rho_2:\Gamma_{\alpha} \rightarrow P^+(4).$$
Moreover, since $\rho_1\in C^{1, 0}(\Gamma_{\alpha}, P_{(0)}^+(4))$,
$$\rho_k\in C^{1, 0}(\Gamma_{\alpha}, P_{(2-2k)}^+(4)) \quad \hbox{ for } k\geq 1.$$
Also, all $\rho_k$ have  weight zero.
A solution of  equation (4.33) is given as follows:
$$\rho_2(F_{\alpha}^1) = t^{-2}.$$
Next,
$$d\rho_3  = - {1\over 2}\sum_{i + j = 3}\lbrack\lbrack\rho_i,\rho_j\rbrack\rbrack.$$
Since $\lbrack\lbrack\rho_1,\rho_2\rbrack\rbrack = 0$,  $\rho_3 = 0$. Also,
$$d\rho_4 = - {1\over 2}\lbrack\lbrack\rho_2,\rho_2\rbrack\rbrack.$$
Clearly, $\lbrack\lbrack\rho_2,\rho_2\rbrack\rbrack = 0$, hence $\rho_4 = 0$.
Thus
$\tilde{\rho}_{\alpha, \beta} = \rho_{\alpha} + \beta\rho_1 + \beta^2\rho_2$
is a formal deformation of the embedding $\rho_{\alpha}$.
$$\eqno\Q$$

\noindent
{\bf Remark 4.4.} Note that there exists an embedding of
$\Gamma (2, -1-{\alpha}, {\alpha}-1)$ into
the Lie superalgebra $P_{\hbox{h}}^+(4)$   of differential operators on $S^{1|2}$,
which contracts to $P^+(4)$.

Let $P_{\hbox{h}}^+$
be an associative algebra having the same vector space as
$P^+$, where $\hbox{h}\in  (0, 1]$ and
the multiplication  is given as follows:
$$A(t, \tau)\circ_{\hbox {h}}B(t, \tau) = \sum_{n\geq 0 }{{\hbox{h}}^n\over {n!}}
\partial_{\tau}^nA(t, \tau)\partial_t^nB(t, \tau).$$
The Lie algebra structure on the vector space $P_{\hbox{h}}^+$ is given by
$$[A, B]_{\hbox{h}} = {1\over {\hbox{h}}}(A\circ_{\hbox{h}}B - B\circ_{\hbox{h}}A),$$
so that
$$\hbox{lim}_{\hbox{h}\rightarrow 0}[A, B]_{\hbox{h}} = \lbrace A, B\rbrace.\eqno (4.34)$$
Let $\Lambda_{\hbox{h}}(2N)$ be an associative superalgebra with generators
$\xi_1, \ldots, \xi_N, \eta_1, \ldots, \eta_N$ and relations
$$\xi_i\xi_j = - \xi_j\xi_i,\quad \eta_i\eta_j = -\eta_j\eta_i, \quad \eta_i\xi_j = \hbox{h}\delta_{i, j} - \xi_j\eta_i.$$
Let $P_{\hbox{h}}^+(2N) = P_{\hbox{h}}^+\otimes\Lambda_{\hbox{h}}(2N)$ be a superalgebra with the product given by
$$(A_1\otimes X)(B_1\otimes Y) = (A_1\circ_{\hbox{h}}B_1)\otimes (XY),$$
 where
$A_1, B_1 \in P_{\hbox{h}}^+$ and $X, Y \in \Lambda_{\hbox{h}}(2N)$.
The Lie bracket of $A = A_1\otimes X$ and $B = B_1\otimes Y$  is
$$[A, B]_{\hbox{h}} = {1\over {\hbox{h}}}(AB - (-1)^{p(X)p(Y)}BA),$$
and (4.34) holds.

\noindent
We proved in [3] that
for each $\hbox{h}\in (0, 1]$ and each $\alpha\in\C$
there exists an  embedding
$$\rho_{\alpha, \hbox{h}}: \Gamma (2, -1 - \alpha, \alpha - 1)
\rightarrow P_{\hbox{h}}^+(4).$$
$\Gamma_{\alpha, \hbox{h}} = \rho_{\alpha, \hbox{h}}(\Gamma (2, -1 - \alpha, \alpha - 1))$
is spanned by the following elements:

\begin{equation*}
\begin{aligned}
&E_{\alpha,\hbox{h}}^1 =  t^2,\quad H_{\alpha,\hbox{h}}^1 = t\tau + {{\alpha + 1}\over 2}\hbox{h},\\
&F_{\alpha,\hbox{h}}^1 = \tau^2 - \alpha(2t^{-2}\xi_1\xi_2\eta_1\eta_2 +
t^{-2}(\xi_1\eta_1 + \xi_2\eta_2)\hbox{h} - t^{-1}\tau \hbox{h}),\\
&E_{\alpha,\hbox{h}}^2 = \xi_1\xi_2,\quad H_{\alpha,\hbox{h}}^2 =  \xi_1\eta_1 + \xi_2\eta_2 - \hbox{h},\quad
F_{\alpha,\hbox{h}}^2 = \eta_1\eta_2,\\
&E_{\alpha,\hbox{h}}^3 = \xi_1\eta_2,\quad H_{\alpha,\hbox{h}}^3 =  \xi_1\eta_1 - \xi_2\eta_2,\quad
F_{\alpha,\hbox{h}}^3 = \xi_2\eta_1,\\
&T_{\alpha,\hbox{h}}^1 =   t\eta_1,\quad
T_{\alpha,\hbox{h}}^2 =  t\eta_2,\quad
T_{\alpha,\hbox{h}}^3 =  t\xi_1,\quad
T_{\alpha,\hbox{h}}^4 =  t\xi_2,\\
&D_{\alpha,\hbox{h}}^1 =  \tau\xi_1 + \alpha t^{-1}\xi_1\xi_2\eta_2,\quad
D_{\alpha,\hbox{h}}^2 = \tau\xi_2 - \alpha t^{-1}\xi_1\xi_2\eta_1,\\
&D_{\alpha,\hbox{h}}^3 = \tau\eta_1 + \alpha t^{-1}\eta_1\eta_2\xi_2,\quad
D_{\alpha,\hbox{h}}^4 = \tau\eta_2 - \alpha t^{-1}\eta_1\eta_2\xi_1,\\
\end{aligned}
\end{equation*}
so that
$\hbox{lim}_{\hbox{h}\rightarrow 0}\Gamma_{\alpha, \hbox{h}} =
\Gamma_{\alpha}\subset P^+(4)$.
Note that in [13] we constructed a different embedding $\Gamma_{\alpha, \hbox{h}}\subset P_{\hbox{h}}(4)$, where we essentially used pseudodifferential symbols.
Using the same methods as in [13],
we obtain the following analogue of Theorem 4.3:

\vskip 0.1in
\noindent
{\bf Theorem 4.5.}
$\hbox{dim}H^1(\Gamma_{\alpha, \hbox{h}}, P^+_{\hbox{h}}(4)) = 1$. The cohomology space is spanned by the class of the 1-cocycle $\bar{\theta}_1$ given as follows:
\begin{equation*}
\begin{aligned}
&\bar{\theta}_1(D_{\alpha,\hbox{h}}^1) =  t^{-1}\xi_1,\quad
\bar{\theta}_1(D_{\alpha, \hbox{h}}^2) =  t^{-1}\xi_2,\\
&\bar{\theta}_1(D_{\alpha, \hbox{h}}^3) =  t^{-1}\eta_1,\quad
\bar{\theta}_1(D_{\alpha, \hbox{h}}^4) =  t^{-1}\eta_2,\\
&\bar{\theta}_1(F_{\alpha, \hbox{h}}^1) = 2t^{-1}\tau +(\alpha - 1)t^{-2}\hbox{h},\quad
\bar{\theta}_1(H_{\alpha, \hbox{h}}^1) = 1.\\
\end{aligned}
\end{equation*}
The infinitesimal deformation defined by $\bar{\theta}_1$
corresponds to the formal deformation of the embedding $\rho_{\alpha, \hbox{h}}$ given as follows:
$\tilde{\rho}_{\alpha, \beta, \hbox{h}} = \rho_{\alpha, \hbox{h}} + \beta\rho_1 + \beta^2\rho_2$, where
$\rho_1 = \bar{\theta}_1$ and $\rho_2$ is defined by
$$\rho_2(F_{\alpha,\hbox{h}}^1) = t^{-2}.$$

\vskip 0.3in
\noindent
{\bf Acknowledgement}
\vskip 0.3in

This work is based on the preprint
of the Max-Planck-Institut f{\"u}r Mathematik
MPIM2008--82.
The author would like to thank the MPIM in Bonn
for the hospitality and support.

\vfil\eject
\noindent
{\bf References}

\vskip 0.2in

\begin {itemize}

\font\red=cmbsy10
\def\~{\hbox{\red\char'0016}}

\item[{[1]}]
V. G. Kac,
Lie superalgebras,
Adv. Math. 26 (1977) 8--96.

\item[{[2]}] M. Scheunert, The Theory of Lie Superalgebras, in:
Lecture Notes in Mathematics 716,  Springer, Berlin, 1979.

\item[{[3]}] E. Poletaeva,
Embedding of the Lie superalgebra $D(2, 1 ; \alpha)$
into the Lie superalgebra of pseudodifferential symbols on $S^{1|2}$,
J. Math. Phys. 48 (2007) 103504, 17 pp.; e-print arXiv:0709.0083.

\item[{[4]}] A. Nijenhuis, R. W. Richardson Jr.,
Deformations of homomorphisms of Lie groups and Lie algebras,
Bull. Amer. Math. Soc. 73 (1967) 175--179.

\item[{[5]}]  R. W. Richardson Jr.,
Deformations of subalgebras of Lie algebras,
J. Diff. Geom. 3 (1969) 289--308.

\item[{[6]}]
V. Ovsienko and C. Roger,
Deforming the Lie algebra of vector fields on $S^1$
inside the Poisson algebra on $\dot{T}^*S^1$,
Comm. Math. Phys. 198 (1998) 97--110.

\item[{[7]}]
V. Ovsienko and C. Roger,
Deforming the Lie algebra of vector fields on $S^1$
inside the Lie algebra of pseudodifferential symbols on $S^1$,
Amer. Math. Soc. Transl.  194 (1999) 211--226.

\item[{[8]}] B. Agrebaoui,  N. Ben Fraj, and S. Omri,
On the cohomology of the Lie superalgebra of contact vector fields on $S^{1|2}$,
J. Nonlinear Math. Phys. 13 (2006)  523--534.

\item[{[9]}]  N. Ben Fraj, S. Omri,
Deforming the Lie superalgebra of contact vector fields on $S^{1|1}$ inside the Lie superalgebra of superpseudodifferential operators on $S^{1|1}$,
J. Nonlinear Math. Phys. 13 (2006)  19--33.

\item[{[10]}]
V. G. Kac,
Classification of supersymmetries, in:
Proceedings of the International Congress of Mathematicians, Beijing, 2002 (Higher Education Press, Beijing, 2002), Vol. I, pp. 319--344.

\item[{[11]}] E. Poletaeva,
A spinor-like representation of the contact
superconformal algebra  $K'(4)$,
J. Math. Phys.  42 (2001) 526--540; e-print arXiv:hep-th/0011100
and references therein.

\item[{[12]}] D. B. Fuchs,
Cohomology of Infinite-Dimensional Lie Algebras, Consultants Bureau, New York, 1987.

\item[{[13]}]
E. Poletaeva,
Deforming the Lie superalgebra $D(2, 1 ; \alpha)$
inside the superconformal algebra $K'(4)$.
Journal of Mathematical Sciences  161, (2009) 130--142.

\end{itemize}

\end{document}